\documentclass[11pt]{amsart}
\usepackage{amsmath,amsfonts,latexsym,amssymb,amscd, graphicx,pictexwd}

\renewcommand{\P}{{\mathbb P}}
\newcommand{\E}{{\mathbb E}}
\newcommand{\R}{{\mathbb R}}
\newcommand{\RR}{{\mathbb R}}

\newcommand{\ZZ}{{\mathbb Z}}
\newcommand{\PP}{{\mathbf P}}

\newcommand{\0}{{\mathbf 0}}

\newcommand{\cc}{{\mathbf c}}
\newcommand{\z}{{\mathbf z}}

\newcommand{\cal}{\mathcal}

\newcommand{\calN}{{\mathcal N}}

\newcommand{\eps}{{\varepsilon }}

\newcommand{\Exp}{{\rm Exp}}

\newtheorem{lemm}{Lemma}[section]

\numberwithin{equation}{section}

\def\x{\mathbf{x}}
\def\t{\mathbf{t}}
\def\v{\mathbf{v}}
\def\y{\mathbf{y}}
\def\D{\mathcal{D}}
\def\bP{\mathbf{P}}
\def\B{\mathbf{B}}
\def\one{{\mathbf 1}}

\newtheorem{theo}{Theorem}

\newtheorem{rem}{Remark}[section]
\newenvironment{dem}{\vskip 2mm\noindent {\it Proof} :}
                    {\hfill $\square$ \vskip 2mm \noindent}

\begin{document}

\title[Busemann functions and the speed of a second class particle]{Busemann functions and the speed of a\\second class particle in the rarefaction fan}

\author{Eric Cator}
\address{Delft University of Technology\\
Mekelweg 4, 2628 CD Delft, The Netherlands}
\email{E.A.Cator@tudelft.nl}

\author{Leandro P. R. Pimentel}
\address{Institute of Mathematics\\
Federal University of Rio de Janeiro}
\email{leandro@im.ufrj.br}
\thanks{Leandro P. R. Pimentel was supported by grant numbers 613.000.605 and 040.11.146 from the Netherlands Organisation for Scientific Research (NWO)}


\begin{abstract}
In this paper we will show how the results found in \cite{CaPi}, about the Busemann functions in last-passage percolation, can be used to calculate the asymptotic distribution of the speed of a single second class particle starting from an arbitrary deterministic configuration which has a rarefaction fan, in either the totally asymetric exclusion process, or the Hammersley interacting particle process. The method will be to use the well known last-passage percolation description of the exclusion process and of the Hammersley process, and then the well known connection between second class particles and competition interfaces.
\end{abstract}

\maketitle

\section{Introduction}

The macroscopic behavior of the density profile of the totally asymmetric nearest neighbor exclusion process (TASEP) is governed by the Burgers equation, which corresponds to the ``shape theorem'' in last-passage percolation \cite{Ro}. The second class particles, that follow roughly the behavior of a perturbation of the system, are asymptotically governed by the characteristics of the Burgers equation. When there is only one characteristic, the second class particle follows it; when there are infinitely many (rarefaction fan), the particle chooses one of them at random to follow. Until recently, there was only one method, introduced by Kipnis and Ferrari \cite{FeKi}, to compute the distribution of the asymptotic speed when the initial distribution is a product measure with densities $p\in[0,1)$ and $p'\in(p,1]$, to the right and left of the origin, respectively. They proved that the rescaled position of the second class particle converges in distribution, as time goes to infinity, to a random variable uniformly distributed in the interval $[1-2p',1-2p]$. Their method is based on the hydrodynamic description of the macroscopic behavior and it strongly relies on the product structure of the initial profile.

Later, following initial results by Guiol and Mountford \cite{GuMo} and Ferrari and Pimentel \cite{FePi}, Ferrari, Martin and Pimentel in \cite{FeMaPi} proved that for initial configurations exhibiting a rarefaction fan, the second class particle almost surely have an asymptotic speed, but no new results were found about the distribution of this speed.

Recently, in a paper by Amir, Angel and Valk\'o \cite{AmAnVa}, a new method has been introduced to study the behavior of a multiclass TASEP, which allows for the calculation of joint distributions of multiclass particles under a specific initial configuration. However, as far as we can see, this does not shed new light on how to calculate the asymptotic distribution of a single second class particle in an arbitrary initial configuration exhibiting a rarefaction fan.

We will also consider the Hammersley interacting particle process \cite{AlDi}, where the situation concerning the second class particle in a rarefaction fan is very similar to TASEP: almost sure existence of the asymptotic speed is proved by Coletti and Pimentel in \cite{CoPi} (they even prove this for more general objects than second class particles), but the distribution of the second class particle in a rarefaction fan was only determined for a very specific family of initial conditions, involving Poisson processes (\cite{CaDo, CoPi}).

In this paper we will show how the results found in \cite{CaPi} about the Busemann functions in last-passage percolation (LPP) can be used to calculate the asymptotic law of the speed of a single second class particle starting from an arbitrary deterministic configuration which has a rarefaction fan, in either the totally asymetric exclusion process (TASEP), or the Hammersley process. The method will be to use the well known LPP description of the TASEP and of the Hammersley process, and then the well known connection between second class particles and competition interfaces (see \cite{FePi} for the TASEP case, and \cite{CaPi} for the Hammersley case) to show that the limit law can be expressed in terms of the supremum of two independent random walks in an environment that depends upon the initial configuration. This method allows us to recover the classical results, as well as to get explicit formulas for the law of the asymptotic speed starting from periodic profiles. Busemann functions in First Passage Percolation were first introduced by Newman \cite{Ne}, and these were later extended to LPP by  W\"uthrich in \cite{Wu}, by Ferrari and Pimentel in \cite{FePi} and by the authors in \cite{CaPi}.

To give some intuition for the main idea in the paper, we present a heuristic outline of it in the lattice LPP competition model. Suppose that to each $\x\in\ZZ^2$ we attach an exponential random variable $X_\x$ of intensity one. For $\x\leq\y$ (coordinate-wise) let $L^\ell(\x,\y)$ denote the maximum, over all up-right nearest-neighbor paths $\varpi$ connecting $\x$ to $\y$, of the sum of $X_\z$ along $\varpi$, leaving out $\x$ (see \eqref{last}). We adopt the notation in \cite{CaPi} and, to indicate that it refers to the lattice model, we use the superscript $\ell$. Now choose two non-ordered points $\x=(x_1,x_2)$ and $\y=(y_1,y_2)$ with $y_1<x_1$ and $y_2>x_2$. We will assign each point $\z\in \R^2$ with $\x,\y\leq \z$ to either $\x$ or $\y$ in the following way: $\x$ ``conquers'' $\z$ if $L^\ell(\x,\z)\geq L^\ell(\y,\z)$, and otherwise $\y$ conquers $\z$. The point $\z$ gets the color blue if it is conquered by $\x$ and otherwise it gets the color red. This example corresponds to the two-corner profile discussed in Section \ref{sec:2corn}. The {\em competition interface} is the up right path that separates the blue and the red regions. It has been shown in \cite{FePi} that a competition interface has an (a.s.) asymptotic angle $\Theta$, which will be random in this case (rarefaction regime). To check whether $\Theta\leq\alpha$, for some fixed $\alpha\in(0,\pi/2)$, one can consider $\z_\alpha$ far away on the line with angle $\alpha$, and compare the length $L^\ell(\x,\z_\alpha)$ with $L^\ell(\y,\z_\alpha)$. It is at this point that the Busemann function comes into play: for $\alpha\in (0,\pi/2)$,
\[ \mbox{ a.s. }\,\,\exists\,\,B^\ell_\alpha(\x,\y):=\lim_{|\z_\alpha|\to +\infty} L^\ell(\y,\z_\alpha) - L^\ell(\x,\z_\alpha)\,.\]
This means that checking whether $L^\ell(\x,\z_\alpha)\leq L^\ell(\y,\z_\alpha)$ (this corresponds to $\Theta\leq \alpha$) eventually, corresponds to checking whether $B^\ell_\alpha(\x,\y)\geq 0$ \footnote{This idea was also used to study the almost sure behavior of competition interfaces \cite{FePi, Pi-2}}. In \cite{CaPi} we derived the distribution of $B^\ell_\alpha$, which is why we are in fact able to calculate the asymptotic distribution of $\Theta$. This argument can be extended for any deterministic initial profile in the rarefaction regime.
\newline

\paragraph{\bf Overview} In Section \ref{results} we state the main results, Theorem \ref{2classTASEP}, Theorem \ref{competition-1} and Theorem \ref{thm:2ndclasspart} . In Section \ref{proofs} we introduce the equilibrium description of the Busemann functions and establish its connection with second class particles and competition interfaces to prove all the theorems.  Section \ref{corollaries} and Section \ref{corollariesH} are devoted to compute the law of the asymptotic speed, for some initial configurations in the rarefaction regime, in the TASEP and in the Hammersley models, respectively.

\section{Main Results}\label{results}

In this section we state the general results about the asymptotic distribution of the second class particle in a rarefaction fan in the TASEP case and the Hammersley case. The TASEP result is formulated in two ways: the first using the particle description, the second using the last passage description. In Section \ref{proofs} we will link the two results. Applications of these general results to particular initial conditions can be found in Section \ref{corollaries} and Section \ref{corollariesH}.

\subsection{Second-class particles in the exclusion process}
The one-dimensional nearest neighbor totally asymmetric simple exclusion process is a
Markov process $(\eta_t\,,\, t\geq 0)$ in the state space $\{0,1\}^\ZZ$. In our notation,
$\eta_t(x)=1$ indicates that there is a particle at site $x$ at time $t$; only one
particle is allowed at each site. At rate one, if there is a particle
at site $x\in\ZZ$, it attempts to jump to $x+1$; if there is no
particle in $x+1$ the jump occurs, otherwise nothing happens. To
construct a realization of this process, one considers
independent one dimensional Poisson processes $\calN=(N_x(\cdot), x\in
\ZZ)$ of intensity 1. The process
$(\eta_t\,,\, t\ge 0)$ can be constructed as a deterministic function
of the initial configuration $\eta$ and the Poisson processes $\calN$ as
follows: if $s$ is a Poisson epoch of $N_x$ and there is a particle at
$x$ and no particle at $x+1$ in the configuration $\eta_{s-}$, then at
time~$s$ the new configuration is obtained by making the particle jump
from $x$ to $x+1$. Let $\Phi$ be the function that takes $\eta$ and
$\calN$ to $(\eta_t\,,\, t\ge 0)$. Let $\eta$ and $\eta'$ be two
arbitrary configurations. The \emph{basic} coupling between two
exclusion processes with initial configurations $\eta$ and $\eta'$
respectively is the joint realization
$(\Phi(\eta,\calN),\Phi(\eta',\calN))=((\eta_t,\eta'_t),\, t\ge 0)$
obtained by using the same Poisson epochs for the two different
initial conditions. Liggett (1985, 1999) are the default references
for the exclusion process.

Let $\eta$ be an initial configuration and let $\eta'$ be the configuration that differs from $\eta$ only at site $0$. With the basic coupling, the configurations at time $t>0$ differ only at one site $X^\eta(t)$ defined by
$$X^\eta(t) := \sum_{x\in\ZZ} x \one\{\eta_t(x)\neq\eta'_t(x)\}\,$$
($\one$ denotes the indicator function). The process $(X^\eta(t),\, t\ge 0)$ is the trajectory of a ``second class particle''. This particle jumps one unit to the right at rate one if there is no $\eta$ particle in its right nearest neighbor and it jumps one unit to the left at rate one if there is an $\eta$ particle in its left
nearest neighbor site, interchanging positions with it.

To state our results, we need some more notation. For each $\eta$ we associate a new configuration $\bar\eta$ as follows:
\begin{equation}\label{etabar}
 \bar\eta(k):=\,\left\{\begin{array}{ll} \eta(k) & \mbox{for } k< 0\\
 0&\mbox{for } k=0\\
 1&\mbox{for }k=1\\
 \eta(k-1) & \mbox{for }k> 1.\end{array}\right.
\end{equation}
The hole-particle pair at the origin will be important to formulate the LPP competition model and to establish the connection between second class particles and competition interfaces \cite{FePi}. Note that $\eta(0)$ is not relevant for $\bar\eta$; this corresponds to the idea of having a second class particle in $0$. For $\rho\in(0,1)$ let $\{\Exp_k(\rho)\,:\,k\in\ZZ\}$ and $\{\Exp'_k(1-\rho)\,:\,k\in\ZZ\}$ be independent collections of i.i.d. exponential random variables of intensity $\rho$ and $1-\rho$, respectively. These collections are also assumed to be independent of $\eta$, whenever $\eta$ is random. For $k> 0$ let
\begin{equation}
 Z^{\rho,\eta}_k:=\,\left\{\begin{array}{ll} -\Exp_k(\rho) & \mbox{if } \bar\eta(k)= 0\\
 \Exp'_k(1-\rho) & \mbox{if }\bar\eta(k)=1\,,\end{array}\right.
\end{equation}
and for $k\leq0$ let
\begin{equation}
Z^{\rho,\eta}_k:=\,\left\{\begin{array}{ll} \Exp_k(\rho) & \mbox{if } \bar\eta(k)= 0\\
 -\Exp'_k(1-\rho) & \mbox{if }\bar\eta(k)=1\,.\end{array}\right.
\end{equation}
Define the random walks
$$S_{n,+}^{\rho,\eta}:=\sum_{k=1}^n Z_k^{\rho,\eta}\mbox { and }S_{n,-}^{\rho,\eta}:=\sum_{k=0}^{-n} Z_k^{\rho,\eta}\,.$$
We assume that the particle system is in the rarefaction regime, i.e.:
\begin{equation}\label{rarefaction}
 p_\eta:=\limsup_{n\to \infty} \frac1{n}\sum_{k=1}^{n}\eta(k) < \liminf_{n\to \infty} \frac1{n}\sum_{k=-n}^{-1} \eta(k)=:p'_\eta\,.
 \end{equation}
\begin{theo}\label{2classTASEP}
For each $u\in [1-2p'_\eta,1-2p_\eta]$ let
$$\rho_u=\frac{1+u}{2}\,.$$
If \eqref{rarefaction} holds then
$$\frac{X^\eta(t)}{t}\stackrel{\D}{\longrightarrow}U^\eta\mbox{ as }t\to\infty\,,$$
where $U^\eta$ is a random variable with support $[1-2p'_\eta,1-2p_\eta]$ and for $u\in (1-2p'_\eta,1-2p_\eta)$
\begin{equation}\label{comp-dist1}
\P\left(U^\eta \geq u\right)= \P\left(\sup_{n\geq 0}S_{n,-}^{\rho_u,\eta}\,\geq\, \sup_{n\geq 1}S_{n,+}^{\rho_u,\eta}\right)\,.
\end{equation}
\end{theo}
To give some intuition for this result, note that the asymptotic drift of the random walk $S_{n,-}^{\rho_u,\eta}$ is lesser than or equal to $-p'_\eta/(1-\rho_u)+(1-p'_\eta)/\rho_u$, which is strictly negative if $u>1-2p'_\eta$ and it is a decreasing function of $u$. Also, the asymptotic drift of the random walk $S_{n,+}^{\rho_u,\eta}$ is lesser than or equal to $-(1-p_\eta)/\rho_u + p_\eta/(1-\rho_u)$, which is strictly negative if $u<1-2p_\eta$ and it is an increasing function of $u$. Note that when $u=1-2p'_\eta$ or $u=1-2p_\eta$, \eqref{comp-dist1} may not hold in general! An example of this can be found in Section \ref{proofs} (see Remark \ref{rem1}). The left-hand side of \eqref{comp-dist1} is still determined at the boundary by left-continuity.

\subsection{Competition interfaces in the lattice last-passage percolation model}
Let $X:=\{X_\z\}_{\x\in \ZZ^2}$ be a collection of i.i.d. random variables with an exponential distribution of para\-meter one. For $\x,\y\in\ZZ^2$, with $\x\leq \y$ (coordinate-wise), let
\begin{equation}\label{last}
L^\ell(\x,\y):=\max_{\varpi:\x\nearrow\y}\sum_{\z\in\varpi}X_\z\,,
\end{equation}
where $\varpi:\x\nearrow\y$ denotes an increasing (or up-right) path connecting $\x$ to $\y$, not including $\x$. The random variable $L^\ell$ is called the last-passage (percolation) time from $\x$ to $\y$.

Using the well known connection between the totally asymmetric exclusion process (TASEP) and the last-passage percolation model, introduced by Rost \cite{Ro}, we identify the initial configuration $\bar\eta$ of particles and holes, that includes a hole at $0$ and a particle at $1$, with an initial (growth) {\em profile} $\sigma$ in $\ZZ^2$, in the following way: we define the profile $\sigma=\sigma_\eta$ (as a function of $\eta$), which is a down-right bi-infinite sequence of points in $\ZZ^2$ containing the origin, such that
\[ \sigma(0)=(0,0)\ \ \mbox{and}\ \ \sigma(k) - \sigma(k-1) = \left\{ \begin{array}{ll}
 (1,0) & \mbox{if } \bar\eta(k) = 0,\\
 (0,-1) & \mbox{if } \bar\eta(k) = 1.
 \end{array}
\right.\]

Given such an initial profile $\sigma$, we can define a competition interface, following \cite{FePi} and \cite{FeMaPi}. The competition for territory goes as follows: choose $\t\in \RR^2$ such that $\{\x\ :\ \x\leq \t\}\cap \sigma \neq \emptyset$. Define
\[L^\ell(\sigma^-,\t):= \sup \{ L^\ell(\sigma(k),\t)\ :\ k<0\ \mbox{and}\ \sigma(k)\leq \t\}\]
and
\[L^\ell(\sigma^+,\t):= \sup \{ L^\ell(\sigma(k),\t)\ :\ k>0\ \mbox{and}\ \sigma(k)\leq \t\}.\]
If $\{\x\ :\ \x\leq \t\}\cap \sigma^\pm = \emptyset$, we define $L^\ell(\sigma^\pm,\t)=0$. The point $\t$ gets the color blue if it is further (the longest path is longer) from the negative part of $\sigma$ than from the positive part, so if
\[ L^\ell(\sigma^-,\t) > L^\ell(\sigma^+,\t),\]
and otherwise it gets the color red. Note that because of the hole-particle pair at the origin, $\sigma(0)$ can never be the most ``distant'' point to $\t$. Also, the coloring is well defined for every $\t$ above $\sigma$ if the initial condition $\eta$ has infinitely many particles to the left of the origin, and infinitely many holes to the right; we will always assume this.
The competition interface $\Phi^\eta:=(\phi^\eta_n)_{n\geq 0}$ is the up right path through the lattice points starting at $\0$ that separates the blue and the red region (see Section 2 of \cite{FeMaPi} for an extensive description).

Under \eqref{rarefaction}, the rarefaction condition, let
$$\alpha_\eta:=\arctan\left(\left(\frac{p_\eta}{1-p_\eta}\right)^2\right)\,<\,\arctan\left(\left(\frac{p'_\eta}{1-p'_\eta}\right)^2\right)=:\alpha'_\eta\,.$$
Note that $p_\eta/(1-p_\eta)$ is the (negative) asymptotic slope of $\sigma$ at $+\,\infty$, whereas $p'_\eta/(1-p'_\eta)$ is the slope of $\sigma$ at $-\,\infty$.
\begin{theo}\label{competition-1}
For each $\alpha\in[\alpha_\eta,\alpha'_\eta]$ let
$$\rho_\alpha:=\frac{1}{1+\sqrt{\tan\alpha}}\,.$$
If \eqref{rarefaction} holds then the competition interface $(\phi^\eta_n)_{n\geq 0}$ satisfies
\begin{equation}\label{comp-direc}
\frac{\phi^\eta_n}{|\phi^\eta_n|}\stackrel{\D}{\longrightarrow}(\cos\Theta^\eta,\sin\Theta^\eta)\mbox{ as }n\to\infty\,,
\end{equation}
where $\Theta^\eta$ is a random variable with support $[\alpha_\eta,\alpha'_\eta]$ and for $\alpha \in (\alpha_\eta,\alpha'_\eta)$
\begin{equation}\label{comp-dist2}
\P\left(\Theta^\eta\leq\alpha\right)=  \P\left(\sup_{n\geq 0}S_{n,-}^{\rho_\alpha,\eta}\,\geq\, \sup_{n\geq 1}S_{n,+}^{\rho_\alpha,\eta}\right)\,.
\end{equation}
\end{theo}
Theorem \ref{competition-1} and Theorem \ref{2classTASEP}, although each of separate interest, are of course intimately linked through the connection between TASEP and the lattice LPP, which enables us to directly translate the competition interface in the LPP into the position of the second class particle in TASEP; see Lemma \ref{2class} in Section \ref{proofs} for the details. Our proof will use the LPP description, as explained in the Introduction. Again we have that at the boundary, \eqref{comp-dist2} may fail.

\begin{rem}\label{rem:femapi}
{\rm
Ferrari, Martin and Pimentel in \cite{FeMaPi} show that the competition interface in the rarefaction fan has in fact almost surely an asymptotic direction (which implies that the second class particle in TASEP has an asymptotic speed). They assume an asymptotic inclination of the initial profile $\sigma$ in the positive and in the negative direction (which corresponds to an asymptotic particle density left and right of zero in TASEP). However, it is not hard to adapt their proof only slightly to get the same results under our rarefaction fan condition \eqref{rarefaction}.}
\end{rem}

\subsection{Second class particles in the Hammersley particle system}
 Let $\PP\subseteq\RR^2$ be a two-dimensional Poisson random set of intensity one. The last-passage time (or longest increasing path) between $\z,\v\in\RR^2$, with $\z\leq\v$, is defined by
\begin{equation}\label{Hammer}
L(\z,\v):=\max_{\varpi:\z\nearrow\v}|\varpi|\,,
\end{equation}
where $\varpi:\z\nearrow\v$ denotes an increasing path from $\z$ to $\v$, consisting of points in $\bP$ that are not on the same horizontal nor vertical level as $\z$. Furthermore, $|\varpi|$ denotes the size (the number of Poisson points) of such a path. The Hammersley model can be seen as a Markov process on the space $\calN$ consisting of the set of all counting processes $\nu$ on $\RR$ such that
\begin{equation}\label{eq:leftden}
\liminf_{y\to-\infty}\frac{\nu(y)}{y}>0\,,
\end{equation}
where
\begin{equation}\label{def-1}
 \nu(y):=\,\left\{\begin{array}{ll} \nu((0,y]) & \mbox{for } y\geq 0\\
-\nu((y,0]) & \mbox{for }y <0.\end{array}\right.
\end{equation}
 (We allow $\nu(y)=-\infty$ for all $y< 0$.) For each $\nu\in\calN$, the process
\begin{equation}\label{last-evo}
L_{\nu}(x,t):=\sup_{z\leq x} \left\{ \nu(z) + L((z,0),(x,t))\right\}\ \ \ \ \ \mbox{ for }\,\,x\in\RR\,,\mbox{ and }\,\,t\geq 0\,
\end{equation}
is well defined and the time evolution defined by
\begin{equation}\label{evolution}
M_t^\nu((x,y]):=L_\nu(y,t)-L_\nu(x,t)\,,
\end{equation}
is a Markov process in $\calN$, called the Hammersley interacting particle system \cite{AlDi}. One can think of $M_t$ as a particle system (or process in $\calN$) that starts with the atoms (or particles) of $\nu$, and that moves according to the (informal) rule: for $(x,t)\in\bP$, at time $t\geq 0$ the nearest particle to
the right of $x$ is moved to $x$, with a new particle created at $x$ if no such particle exists.

To introduce a second class particle in the system, we put an atom at $0$ in the process $\nu$ by defining
\[\bar\nu([0,x])=\nu([0,x])+1\,\,\,\mbox{ for }\,\,\,x\geq 0\,.\]
With this new process, and using the same $\bP$, we construct $M^{\bar\nu}_t$, and then define the location of the second class particle $Y^\nu(t)$ as
\[ Y^\nu(t) = \inf\{x\geq 0: M^{\bar\nu}_t(x)=M^\nu_t(x)+1\}\,.\]
We note that $Y^\nu(t)$ is a non-decreasing function of $t$, meaning that the second class particle moves to the right. We use the notation of \cite{CaPi}, where relevant details can also be found, including the following important connection between the longest path description and the second class particle. Let $\nu^+,\nu^-$ be the processes defined by
\begin{equation}\label{+-}
\nu^+(x)=\left\{\begin{array}{ll} \nu(x) & \mbox{for } x\geq 0\\
-\infty & \mbox{for }x <0.\end{array}\right.
\,\,\mbox{and}\,\,
\nu^-(x)=\left\{\begin{array}{ll} 0 & \mbox{for } x\geq 0\\
\nu(x) & \mbox{for }x <0.\end{array}\right.
\end{equation}
Then
\begin{equation}\label{eq:secondclass}
\{ Y^\nu(t)\leq x\} = \{L_{\nu^+}(x,t)\geq L_{\nu^-}(x,t)\}\,.
\end{equation}

The above equation resembles the situation where one should check if the competition interface is to the left or to the right of some angle $\alpha$. This indicates that the speed of a second class particle can be related to the respective Busemann functions in a similar way. We require the following rarefaction assumption (recall \eqref{def-1} and \eqref{+-})
\begin{equation}\label{eq:assumpint}
a_\nu:=\limsup_{x\to \infty} \frac{\nu^+(x)}{x} < \liminf_{y\to -\infty} \frac{\nu^-(y)}{y}=:b_\nu\,.
\end{equation}
For $\rho>0$ let $\bar{\nu}_\rho$ be a one-dimensional Poisson counting process of intensity $\rho>0$, that is independent of $\nu$, whenever $\nu$ is random.

\begin{theo}\label{thm:2ndclasspart}
For each $v\in[b_\nu^{-2},a_\nu^{-2}]$ let
$$\rho_v:=\frac{1}{\sqrt{v}}\,.$$
If \eqref{eq:assumpint} holds then
$$\frac{Y^\nu(t)}{t}\stackrel{\D}{\to}V^\nu\mbox{ as }t\to\infty\,,$$
where $V^\nu$ is a random variable with support $[b_\nu^{-2},a_\nu^{-2}]$ and for $v\in (b_\nu^{-2},a_\nu^{-2})$,
\begin{equation}\label{2class-dist}
\P\left(V^\nu \leq v\right) = \P\left(\sup_{z\geq 0}  \{\nu(z) -\bar\nu_{\rho_v}(z)\} \geq \sup_{z<0} \{\nu(z) -\bar\nu_{\rho_v}(z)\}\right)\,.
\end{equation}
\end{theo}
Condition \eqref{eq:assumpint} ensures that with probability 1, at least one of the two suprema is finite, and hence the right hand side probability is well defined for every $v\in (b_\nu^{-2},a_\nu^{-2})$. Again we have that at the boundary, \eqref{2class-dist} may fail.

\section{Proof of the Theorems}\label{proofs}
We will start with the proof for the Hammersley process, since the ideas are most clearly represented in that case. In the TASEP case, we will use very similar methods.
\subsection{Busemann functions and equilibrium measures in the Hammersley model}
In the Hammersley last-passage percolation model, W\"{u}thrich \cite{Wu} proved that, for $\alpha\in (0,\pi/2)$,
\[ \mbox{ a.s. }\,\,\exists\,\,B_\alpha(\x,\y):=\lim_{|\z_\alpha|\to +\infty} L(\y,\z_\alpha) - L(\x,\z_\alpha)\,.\]
Here, $\z_\alpha$ follows the ray through $(\cos\alpha,\sin\alpha)$. By the invariance under the map $\alpha\to\alpha+\pi$, it is straightforward to see that the analogous result holds for $\alpha\in(\pi,3\pi/2)$. The most important property of the Busemann functions is the following connection with the associated interacting particle process (Theorem 5.3 and Corollary 5.5 of \cite{CaPi}): Let $\alpha\in(\pi,3\pi/2)$  and define
$$\nu_\alpha((x,y]):=B_{\alpha}(\0,(y,0))-B_{\alpha}(\0,(x,0)) = B_\alpha((x,0),(y,0))\,\,\,\mbox{ for }\,\,\,x<y\,.$$
Then $\nu_\alpha$ is the unique spatially ergodic counting process with intensity
$$\E\nu_\alpha(1)=\sqrt{\tan\alpha}\,,$$
that satisfies (time invariance)
\begin{equation}\label{inv}
M_t^{\nu_\alpha} \stackrel{\D}{=}\nu_\alpha\,,\mbox{ for all }\,t\geq 0\,.
\end{equation}
A consequence of this result is that we can compute the law of the Busemann function along the horizontal axis. For $\rho>0$ let $\bar{\nu}_\rho$ be an one-dimensional Poisson counting process of intensity $\rho>0$. Then
\begin{equation}\label{Poisson}
\nu_\alpha\stackrel{\D}{=}\bar{\nu}_{\sqrt{\tan\alpha}}\,.
\end{equation}
We also need the following lemma for the Hammersley process $L_\nu$ (see \eqref{last-evo}). This result, in a completely different formulation, would also follow from Lemma 2.2 in \cite{CaGr0}.
\begin{lemm}\label{lem:LodominateLnu}
Consider the Hammersley process $L_\nu$, and suppose that for $t\geq 0$, $0\leq x\leq y< Y^\nu(t)$. Then
\[ L_\nu(y,t)-L_\nu(x,t)\leq L(\0,(y,t)) - L(\0,(x,t)).\]
\end{lemm}
\begin{dem}
The fact that $x\leq y<Y^\nu(t)$, implies that (see \eqref{eq:secondclass}) there exists $t_x\leq t_y<0$ such that
\[ L_\nu(x,t) = L((t_x,0),(x,t)) + \nu(t_x)\ \ \ \mbox{and}\ \ \ L_\nu(y,t) = L((t_y,0),(y,t)) + \nu(t_y).\]
Define $\z$ as the (an) intersection point of the longest upright path from $(t_y,0)$ to $(y,t)$ and the longest upright path from $\0$ to $(x,t)$. Then we can see that
\begin{eqnarray*}
L_\nu(x,t) &\geq & L((t_y,0),(x,t)) + \nu(t_y)\\
& \geq & L((t_y,0),\z) + L(\z,(x,t)) + \nu(t_y).
\end{eqnarray*}
This implies
\begin{eqnarray*}
L_\nu(y,t) - L_\nu(x,t) & \leq & L((t_y,0),(y,t)) - L((t_y,0),\z) - L(\z,(x,t))\\
& = & L(\z,(y,t)) - L(\z,(x,t))\\
& = & L(\z,(y,t)) - L(\0,(x,t)) + L(\0,\z)\\
& \leq & L(\0,(y,t)) - L(\0,(x,t)).
\end{eqnarray*}
\end{dem}

\subsection{Proof of Theorem \ref{thm:2ndclasspart}}
For each $v,t>0$  let $\alpha=\alpha_v=\arctan(1/v)$, $x=x_t:=t/ \tan\alpha$, $\x_\alpha := (x,x\tan\alpha)$. By (\ref{eq:secondclass}),
\begin{eqnarray}
\nonumber \{ Y^\nu(t)\leq x\} & = & \{ L_{\nu^+}(\x_\alpha)\geq L_{\nu^-}(\x_\alpha)\}\\
\nonumber& = & \Big\{ \sup_{0\leq z\leq x} \left[ \nu^+(z) + L((z,0),\x_\alpha)\right] \geq \sup_{ z< 0} \left[ \nu^-(z) + L((0,z),\x_\alpha)\right]\Big\}\\
\nonumber& = & \Big\{ \sup_{0\leq z\leq x} \left[ \nu^+(z) + L((z,0),\x_\alpha) - L(\0,\x_\alpha)\right] \\
\label{2classBuse}& & \ \ \ \ \ \ \ \geq \sup_{z<0} \left[ \nu^-(z) + L((0,z),\x_\alpha)- L(\0,\x_\alpha)\right]\Big\}\,.
\end{eqnarray}

For any compact set $K\subset \RR^2$
\begin{equation}\label{eq:busconv}
\left(L(\z,\x_\alpha) - L(\0,\x_\alpha)\,;\, \z\in K\right) \stackrel{\D}{\longrightarrow} \left(B_\alpha(\0,\z)\,;\, \z\in K\right)\,\,\mbox{ as }\,\,x\to \infty\,.
\end{equation}
This remark is the core of the proof. It follows from Theorem 2.2 in \cite{CaPi}: we can take $n\geq 1$ big enough, so that $K\subset [-n,n]\times [-n,n]$. The $\alpha$-rays starting at $(-n,n)$ and $(n,-n)$ will coalesce at some point $\cc\in \R^2$. Furthermore, the longest paths to $\x_\alpha$, starting at $(-n,n)$ and $(n,-n)$, will converge to the respective $\alpha$-rays in a bigger bounded square, containing $\cc$. From that time on, for all $\z\in K$, we will have
\[   L(\z,\x_\alpha) - L(\0,\x_\alpha) = B_\alpha(\0,\z).\]
What follows are some technicalities to show that the location of the maximum behaves well.

\begin{lemm}\label{2classcontrol}
If $\sqrt{\tan\alpha}\in\left(a_\nu,b_\nu\right)$ (recall \eqref{eq:assumpint}) then, as $x\to\infty$,
\[ \sup_{0\leq z\leq x} \left[ \nu^+(z) + L((z,0),\x_{\alpha}) - L(\0,\x_{\alpha})\right]\stackrel{\D}{\longrightarrow} \sup_{z\geq 0} \left[ \nu^+(z) - \nu^{+}_{\alpha}(z)\right]\,,\]
and
\[ \sup_{z< 0} \left[ \nu^-(z) + L((z,0),\x_{\alpha}) - L(\0,\x_{\alpha})\right]\stackrel{\D}{\longrightarrow} \sup_{z< 0} \left[ \nu^-(z) - \nu^{-}_{\alpha}(z)\right].\]
\end{lemm}
\begin{dem}
We will prove that $a_\nu < \sqrt{\tan\alpha}$ implies the first statement. The proof of the second statement follows exactly the same reasoning.

Pick $a<\sqrt{\tan(\alpha)}$ and $C>0$ such that $\nu^+(x)\leq ax + C$ for all $x\geq 0$. For fixed $x$ we know that
\begin{equation}\label{eq:auxdefLnu}
\left(L(\0,\x_\alpha )-L((z,0),\x_\alpha)\,;\,0\leq z\leq x\right) \stackrel{\D}{=} \left(L(\0,\x_\alpha) - L(\0,\x_\alpha-(z,0))\,;\,0\leq z\leq x\right).
\end{equation}
Choose $\alpha'$ and $a'$ such that $a<a'=\sqrt{\tan\alpha'}<\sqrt{\tan\alpha}$. Start a stationary Hammersley system with a Poisson process of intensity $a'$, $\bar\nu_{a'}$, using the same Poisson process in the upper half-plane as $L$. Define as before
\[L_{\bar\nu_{a'}}(x,t) = \sup_{z\leq x}\left\{L((z,0),(x,t)) + \bar\nu_{a'}(z)\right\}.\]
We know that the characteristic direction of this process is $\alpha'$, and from \cite{CaGr} we conclude that the second class particle of this process will be to the right of $\x_\alpha$ with probability greater than $1-O(x^{-3})$; call this event $E_x$. Fix $0\leq z\leq x$. Because of Lemma \ref{lem:LodominateLnu}, under the event $E_x$, we have
\[ L(\0,\x_\alpha) - L(\0,\x_\alpha-(z,0)) \geq L_{\bar\nu_{a'}}(\x_\alpha) - L_{\bar\nu_{a'}}(\x_\alpha - (z,0)).\]
As a process in $z$, the right hand side is just a Poisson process with intensity $a'$. This process will be above the line $az+C$ eventually, and therefore be larger than $\nu^+(z)$. Using \eqref{eq:auxdefLnu}, this means that for each $\eps>0$, we can find $R>0$ such that for all $x$ large enough,
\[ \P\left( \sup_{0\leq z\leq x} \left[ \nu^+(z) + L((z,0),\x_\alpha) - L(\0,\x_\alpha)\right]  >  \sup_{0\leq z\leq R} \left[ \nu^+(z) + L((z,0),\x_\alpha) - L(\0,\x_\alpha)\right]\right) < \eps.\]
Now we can use \eqref{eq:busconv} to finish the proof of the first statement of Lemma \ref{2classcontrol}.
\end{dem}

Notice that, since $\nu_\alpha$ is a function of the underlying two dimensional Poisson random set $\PP$, it is independent of $\nu$. Together with \eqref{Poisson}, \eqref{2classBuse} and Lemma \ref{2classcontrol}, this proves the theorem for $\tan\alpha\in(a_\nu,b_\nu)$.  Now, if $b_\nu <\sqrt{\tan\alpha}$ then a.s. $\sup_{y< 0} \nu^-(y)-\nu^{-}_{\alpha}(y) = \infty$. By the same reasoning as in the  proof of the last lemma,
$$\sup_{z\leq 0} \left[ \nu^-(z) + L((z,0),\x_{\alpha}) - L(\0,\x_{\alpha})\right]\stackrel{\D}{\longrightarrow} \infty\,\,,\mbox { as }\,\,x\to\infty\,.$$
Hence, both functions in \eqref{2class-dist} are zero. The case where $\sqrt{\tan\alpha} < a_\nu$ follows from a similar argument.

\begin{rem}\label{rem1}
{\rm Note that the limiting distribution of $t^{-1}Y^\nu(t)$, for a deterministic $\nu$, may only have an atom at $\alpha$ such that $\sqrt{\tan\alpha} = a_\nu$ or $\sqrt{\tan\alpha} = b_\nu$. If this is the case, then the function \eqref{2class-dist} is not continuous at these choices of $\alpha$, but it is continuous everywhere else. An example would be if $\nu^+(y) = \lfloor y\rfloor - \lfloor {\lfloor y\rfloor}^{2/3}\rfloor$. In this case, when we take $\alpha=\pi/4$, which corresponds to $a_\nu=\sqrt{\tan\alpha}$, and we make sure that $a_\nu<b_\nu$, we see that
\[ \P\left(\sup_{y\geq 0} \left[\nu^+(y) - \nu_\alpha(y)\right] \geq \sup_{y<0}\left[ \nu^-(y) - \nu_\alpha(y)\right]\right) <1,\]
whereas $\P(\lim_{t\to \infty} Y^{\nu}(t)/t\leq 1/\tan\alpha)=1$. This means that $\P(\lim_{t\to \infty}Y^{\nu}(t)/t= 1)>0$.}
\end{rem}

\subsection{Busemann functions and equilibrium measures in the TASEP-LPP context}
For the lattice last-passage percolation model with exponential weights, the existence of the Busemann functions was proved in \cite{FePi}: for $\alpha\in (0,\pi/2)$,
\[ \mbox{ a.s. }\,\,\exists\,\,B^\ell_\alpha(\x,\y):=\lim_{|\z_\alpha|\to +\infty} L^\ell(\y,\z_\alpha) - L^\ell(\x,\z_\alpha)\,.\]
By the invariance under the map $\alpha\to\alpha+\pi$, it is straightforward to see that the analog result holds for $\alpha\in(\pi,3\pi/2)$. In this set up  one can also define last-passage times $L_{\nu^\ell}$ similar to \eqref{last-evo}, with respect to an initial configuration $\nu^\ell$ defined on $\ZZ$, and an evolution $M$ similar to \eqref{evolution}. This evolution will be a Markov process and
$$\nu^\ell_\alpha((x,y]):=B^\ell_{\alpha}(\0,(y,0))-B^\ell_{\alpha}(\0,(x,0))\,\,\,\mbox{ for }\,\,\,x<y\,,$$
will be the unique spatially ergodic process with intensity
$$\E\nu^\ell_\alpha(1)=\frac{1}{1+\sqrt{\tan\alpha}}\,(=\rho_\alpha)\,,$$
that is time invariant (Theorem 8.1 of \cite{CaPi}). This allow us to compute the distribution of $B^{\ell}$ as follows:

\begin{lemm}\label{BuseDist}
The Busemann function is additive,
$$\B^\ell_\alpha(\x,\y)=\B^\ell_\alpha(\x,\z)+\B^\ell_\alpha(\z,\y)\,,$$
and anti-symmetric,
$$\B^\ell_\alpha(\x,\y)=-\B^\ell_\alpha(\y,\x)\,.$$
Suppose $\alpha\in(0,\pi/2)$. Define $Z_k = B^\ell_\alpha(\sigma(k-1),\sigma(k))$. Then all $Z_k$'s are independent and
\begin{equation}\label{lem:profile}
Z_k  \stackrel{\D}{=}\left\{ \begin{array}{ll}
{\rm -\Exp}\left\{\rho_\alpha\right\} & \mbox{if } \bar{\eta}(k) = 0\\
{\rm \Exp}\left\{1 - \rho_\alpha\right\} & \mbox{if }\bar\eta(k) = 1\,,
\end{array} \right.
\end{equation}
where $\Exp\left\{\rho\right\}$ denotes a exponential random variable of intensity $\rho$.
\end{lemm}

\begin{dem}
The first two statements are general properties of the Busemann function and follow from Proposition 3.1 of \cite{CaPi}. Property \eqref{lem:profile} follows from Theorem 8.1 of \cite{CaPi} (see also \cite{BaCaSe}) together with anti-symmetry.
\end{dem}

\subsection{Proof of Theorem \ref{competition-1}}
For each $t>0$ let $\varphi^\eta(t)$ be the intersection of the competition interface $\Phi^\eta=(\phi_n)_{n\geq 1}$ with the line $\{(x,t)\,:\,x\geq 0\}$. For each $\alpha\in(0,\pi/2)$ let $x=x_t:=t/ \tan\alpha$, $\x_\alpha := (x,x\tan\alpha)$. Then
\begin{eqnarray}
\nonumber \{ \varphi^\eta(t)\leq x\} & = & \{ L^\ell(\sigma^+,\x_\alpha)\geq L^\ell(\sigma^-,\x_\alpha)\}\\
\nonumber& = & \Big\{ \sup_{k> 0, \sigma(k)\leq \x_\alpha}L^\ell(\sigma(k),\x_\alpha)\geq \sup_{ k< 0,\sigma(k)\leq \x_\alpha}L^\ell(\sigma(k),\x_\alpha)\Big\}\\
\nonumber& = & \Big\{ \sup_{k> 0, \sigma(k)\leq \x_\alpha} \left[L^\ell(\sigma(k),\x_\alpha) - L^\ell(\0,\x_\alpha)\right] \\
\nonumber& & \ \ \ \ \ \ \ \geq \sup_{ k< 0,\sigma(k)\leq \x_\alpha} \left[L^\ell(\sigma(k),\x_\alpha)- L^\ell(\0,\x_\alpha)\right]\Big\}\,.
\end{eqnarray}
We now need to control the location of the maximum on the profile. For this we use the following lemma:
\begin{lemm}\label{lem:prop31}
Suppose $\alpha<\alpha'_\eta$. Almost surely, there exists $M>0$ such that for all $x>0$,
\[ \sup_{k<0}\ L^\ell(\sigma(k),\x_\alpha) = \sup_{-M\leq k<0}\ L^\ell(\sigma(k),\x_\alpha).\]
Similarly, suppose $\alpha>\alpha_\eta$. Almost surely, there exists $M>0$ such that for all $x>0$,
\[\sup_{k>0}\ L^\ell(\sigma(k),\x_\alpha) = \sup_{0<k\leq M}\ L^\ell(\sigma(k),\x_\alpha).\]
\end{lemm}
\begin{dem}
This result follows more or less directly from Proposition 3.1 of \cite{FeMaPi}. That result is stated for profiles that have an asymptotic direction, whereas we only have \eqref{rarefaction}. However, the proof in the more general case is exactly the same. The proposition states (translated to our notation), that if $\alpha$ and $\alpha'$ are chosen such that $\tan\alpha > \tan\alpha' > \tan\alpha_\eta$, and $\{\z_i\}$ is a sequence of points going to infinity with asymptotic direction $\alpha'$, then there are only finitely many $\z_i$'s that belong to a longest path from $\sigma_+$ to any of the $\x_\alpha$'s. If we take $M_+>0$ such that the $x$-coordinates of those finitely many $\z_i$'s are smaller than the $x$-coordinate of $\sigma(M_+)$, then it is clear that all the longest paths from $\sigma_+$ to any of the $\x_\alpha$'s must start before $\sigma(M_+)$. We apply Proposition 3.1 of \cite{FeMaPi} again to the negative part of $\sigma$, obtaining $M_->0$. Now put $M=\max\{M_-,M_+\}$.
\end{dem}
The proof now follows as in the Hammersley case, using the convergence to the Busemann function on a compact set. Therefore,
\begin{equation}\label{angle}
\P(\Theta^\eta\leq \alpha)=\P\left( \sup_{k<0}B^\ell_\alpha(\0,\sigma(k))\geq  \sup_{k>0}B^\ell_\alpha(\0,\sigma(k))\right)\,,
\end{equation}
and Theorem \ref{competition-1} is a consequence of Lemma \ref{BuseDist} together with \eqref{angle}, for $\alpha\in (\alpha_\eta,\alpha'_\eta)$.

We need an argument for $\alpha<\alpha_\eta$ (the argument for $\alpha>\alpha'_\eta$ will be analogous). When considering $L(\sigma_-,\x_\alpha)$, Lemma \ref{lem:prop31} tells us that almost surely, there exists $M>0$ such that
\[ L(\sigma_-,\x_\alpha)\leq L((-M,0),\x_\alpha).\]
Choose $0<p<p_\eta$ such that 
\[ \rho_\alpha := \frac1{1+\sqrt{\tan\alpha}} >1 - p.\]
This is possible, since $\alpha<\alpha_\eta$ implies that $\rho_\alpha > 1 - p_\eta$. Define a sequence $\tilde\sigma_+(k)=(k,\lfloor -\lambda_p k\rfloor)$, with $\lambda_p= p/(1-p)$. Because of the rarefaction assumption \eqref{rarefaction}, we know that there are infinitely many points on $\sigma_+$ below $\tilde\sigma_+$. We first show that with probability 1, there exists $k_0\geq 1$ such that for $k\geq k_0$,
\begin{equation}\label{eq:Busepos}
B_\alpha((-M,0),\tilde\sigma_+(k))>0.
\end{equation} 
This follows from the fact that
\begin{eqnarray*}
B_\alpha((-M,0),(k,\lfloor -\lambda_p k\rfloor)) & = & B_\alpha((-M,0),(k,0)) + B_\alpha((k,0),(k,\lfloor -\lambda_p k\rfloor)).
\end{eqnarray*}
From Lemma \ref{BuseDist}, we know that $B_\alpha((-M,0),(k,0))$ is distributed as the sum of $k+M$ independent exponentials with expectation $-1/\rho_\alpha$, whereas $B_\alpha((k,0),(k,\lfloor -\lambda_p k\rfloor)$ is distributed as the sum of $\lceil\lambda_p k\rceil$ independent exponentials with expectation $1/(1-\rho_\alpha)$. Now \eqref{eq:Busepos} follows for $k$ big enough from the observation that
\[ -\frac1{\rho_\alpha} + \frac{\lambda_p}{1-\rho_\alpha} > -\frac{1}{1-p} + \frac{\lambda_p}{p} =0.\]
Now choose $\z\in \sigma_+$ such that $\z\leq \tilde\sigma_+(k)$, for some $k\geq k_0$. Then $B_\alpha((-M,0),\z)> 0$. Therefore, for $x$ large enough, we will have that $L(\z,\x_\alpha) > L((-M,0),\x_\alpha)\geq L(\sigma_-,\x_\alpha)$. Since $\z\in \sigma_+$, this means that with probability 1, the angle of the competition interface will eventually be greater than $\alpha$.

\subsection{Proof of Theorem \ref{2classTASEP}} Together with Lemma \ref{2class} below, Theorem \ref{competition-1} implies Theorem \ref{2classTASEP}. This lemma introduces the hole-particle representation of the second class particle to make the link with the competition interface.
\begin{lemm}\label{2class}
For $u\in(-1,1)$ let
$$\alpha_u:=\arctan\left(\frac{1-u}{1+u}\right)^2\,.$$
Then
$$\P\left(U^\eta\geq u\right)=\P\left(\Theta^\eta\leq \alpha_u\right)\,.$$
\end{lemm}
\begin{dem}
We label the particles sequentially from right to
left and the holes from left to right, with the convention that the
particle at site one and the hole at site zero are both labeled $0$.  Let
$P_j(0)$ and $H_j(0)$, $j\in\ZZ$ be the positions of the particles and holes
respectively at time 0 of the initial configuration $\bar\eta$. The position at time $t$ of the $j$th particle
$P_{j}(t)$ and the $i$th hole $H_{i}(t)$ are functions of the variables $L(\sigma,\z)$: at time $L(\sigma,(i,j))$, the $j$th particle and the $i$th hole
interchange positions. Disregarding labels and defining
$\bar\eta_{t}(P_{j}(t))=1,\,\bar\eta_{t}(H_{j}(t))=0,\, j\in\ZZ$, the process $\bar\eta_t$
indeed realizes the exclusion dynamics (see \cite{Ro}). Now, let $\tau_0=0$ and
$$\tau_n:=\max\{L(\sigma^-,\phi_n),L(\sigma^+,\phi_n)\}\,.$$
Define the process $((I^\sigma(t),J^\sigma(t))\,,\,t\geq 0)$ by
$$(I^\sigma(t),J^\sigma(t))=\phi_n\,\mbox{ for }t\in[\tau_n,\tau_{n+1})\,.$$
This process is the hole-particle representation of the second class particle. By Proposition 3 of \cite{FePi},
$$(X^\eta(t)\,,\,t\geq 0) \stackrel{\D}{=}(I^\sigma(t)-J^\sigma(t)\,,\,t\geq 0)\,.$$
By using a similar argument as in the proof of   Proposition 5 of \cite{FePi}, one has that
$$\frac{I^\sigma(t)-J^\sigma(t)}{t} \stackrel{\D}{\to}f(\Theta)\,\mbox{ as }t\to\infty$$
where
$$f(\alpha):=\frac{1-\sqrt{\tan\alpha}}{1+\sqrt{\tan\alpha}}\,,$$
which implies Lemma \ref{2class}.
\end{dem}

\section{Computing the law of the asymptotic speed in the TASEP}\label{corollaries}

\subsection{Two-corner initial profile.}\label{sec:2corn} Choose integers $x,y\geq 1$ and let $\eta=\eta_{x,y}$ be defined as follows:
\[ \eta_{x,y}(k) := \left\{ \begin{array}{ll}
0 & \mbox{if } -(x-1)\leq k<0 \,\mbox{ or } k\geq y\\
1 &  \mbox{if } 0< k \leq y-1\,\mbox{ or } k\leq -x.
\end{array}
\right. \]
For this initial configuration the supremum is attained at deterministic points and it is easy to see that
$$\sup_{n\geq 0}S_{n,-}^{\rho,\eta}\stackrel{\D}{=}\Gamma_{x,\rho}\,\mbox{ and }\,\sup_{n\geq 0}S_{n,+}^{\rho,\eta}\stackrel{\D}{=}\Gamma_{y,1-\rho}\,,$$
where $\Gamma_{x,\rho}$ and $\Gamma_{y,1-\rho}$ are two independent gamma random variables of parameters $(x,\rho)$ and $(y,1-\rho)$, respectively. On the other hand,
$$\P\left(\Gamma_{x,\rho}\geq\Gamma_{y,1-\rho}\right)= \frac{1}{\Gamma(x)\Gamma(y)}\int_0^\infty \Gamma\left(x,\frac{\rho}{1-\rho}z\right)z^{y-1}e^{-z}dz\,,$$
where, for $x,z\in[0,\infty)$,
$$\Gamma(x,z):=\int_z^\infty u^{x-1}e^{-u}du=\frac{\Gamma(x)}{e^{z}}\sum_{j=0}^{x-1}\frac{z^j}{j!}\,\,\mbox{ and }\,\,\Gamma(x):=\Gamma(x,0)=(x-1)!\,.$$
By calculating the integral, one finds that
$$\P\left(\Gamma_{x,\rho}\geq\Gamma_{y,1-\rho}\right)=\frac{1}{(y-1)!}(1-\rho)^y\sum_{j=0}^{x-1}\frac{(j+y-1)!}{j!}\rho^j\,,$$
which implies that
$$\P\left(U\geq u\right)=\frac{1}{(y-1)!}\left(\frac{1-u}{2}\right)^y\,\,\sum_{j=0}^{x-1}\frac{(j+y-1)!}{j!}\left(\frac{1+u}{2}\right)^j\,,$$
and that
$$\P\left(\Theta\leq\alpha\right)=\frac{1}{(y-1)!}\left(\frac{\sqrt{\tan\alpha}}{1+\sqrt{\tan\alpha}}\right)^y\,\,\sum_{j=0}^{x-1}\frac{(j+y-1)!}{j!}\left(\frac{1}{1+\sqrt{\tan\alpha}}\right)^j\,.$$
Compare with the results found in \cite{FeKi,FePi} for $x=y=1$. We also remark that, for $x=L$ and $y=1$, $U\stackrel{\D}{\to} 1\ ({L\to\infty})$, while for $x=1$ and $y=L$, $U\stackrel{\D}{\to} -1\ ({L\to\infty})$. Furthermore, if $x=y=L$, then $U\stackrel{\D}{\to} 0\ ({L\to\infty})$; this follows immediately from the fact that $\E(\Gamma(L,\rho)-\Gamma(L,1-\rho)) = O(L)$ if $\rho\neq 1/2$, whereas the variance is also of order $L$. This answers a question posed in \cite{FeMaPi} (see the end of Section 2.2 there).
\newline

\subsection{Periodic initial configuration.}\label{sec:k+} We consider the following initial configuration $\eta$ for TASEP: choose $k_+,k_-\geq 1$, and suppose to the right of zero we start with $k_+$ holes, then one particle, then $k_+$ holes, etc. So for the configuration $\bar\eta$ with a hole-particle pair at zero we have, for $i\geq 1$,
\[ \bar\eta(i)=\left\{ \begin{array}{ll}
1 & {\rm if}\ i\,{\rm mod}\,(k_++1) = 1,\\
0 & {\rm if}\ i\,{\rm mod}\,(k_++1)\neq 1.
\end{array}\right.\]
To the left of zero, we have a similar, dual pattern: we start with $k_-$ particles, then one hole, then $k_-$ particles, etc. So for $i\leq 0$,
\[ \bar\eta(i)=\left\{ \begin{array}{ll}
0 & {\rm if}\ i\,{\rm mod}\,(k_-+1) = 0,\\
1 & {\rm if}\ i\,{\rm mod}\,(k_-+1)\neq 0.
\end{array}\right.\]
If $\max\{k_+,k_-\}\geq 2$, then we are in the rarefaction regime. The profile $\sigma$ corresponding to $\eta$ starts at $(0,0)$ with a jump to the right, then $k_+$ jumps down, then a jump to the right, etc. Going to the left, it starts with a jump up, then $k_-$ jumps to the left, then a jump up, etc.
\newline

\paragraph{\bf The G/M/1 queue}
The key to solving the distribution of
$$S_+:=\sup_{n\geq 0}S^{\rho,\eta}_{n,+}\,\,\mbox{  and }\,\,S_-:=\sup_{n\leq -1}S^{\rho,\eta}_{n,-}\,$$
is the following connection to a G/M/1 queue: suppose the arrival times $A_i$ of the queue have a distribution
\[ A_i \sim \Exp\{\rho\} + \ldots + \Exp\{\rho\}\ \ \ \ \ ({\rm i.i.d. }\,k_+\ {\rm times})\,.\]
The server time $B_i$ is exponential, $B_i\sim \Exp\{1-\rho\}$. If we define $X_i=B_i-A_i$, and $Z_+\sim \Exp\{1-\rho\}$ independent of the $X_i$'s, then
\[ S_+ \sim Z_+ + \max\left\{0\,,\,\sup_{n\geq 2} \sum_{i=1}^{n-1} X_i \right\}\,.\]
It is a well known fact (see for example Part II of \cite{Co}) that the second term of the right-hand side has the same distribution as the waiting time $W_\infty$ of the stationary G/M/1 queue, which exists if $\E A_i>\E B_i$, or
\[ \frac{1-\rho}{\rho} > \frac1{k_+} \Leftrightarrow 1-\rho > \frac1{1+k_+}\,.\]
Note that $1/(1+k_+)$ equals the fraction of particles to the right of zero. The distribution of $W_\infty$ is explicitly known in terms of the Laplace transform of $A_i$: define
\[ \phi(s) = \E\left(e^{-sA_1}\right) = \left(\frac{\rho}{\rho+s}\right)^{k_+}\,.\]
Now define $\lambda_+=\lambda_+(\rho)$ as the smallest positive zero of
\[ z\mapsto z-\phi((1-\rho)(1-z))\,.\]
Then $\lambda_+$ is the smallest positive solution of
\begin{equation}\label{periodic+}
 \lambda_+(1 - (1-\rho) \lambda_+)^{k_+} = \rho^{k_+}\,.
 \end{equation}
Note that $\lambda_+\uparrow 1$ as $1-\rho\downarrow 1/(1+k_+)$. Given this constant, we know that for all $s\geq 0$,
\[ \P(W_\infty>s) = \lambda_+e^{-(1-\rho)(1-\lambda_+)s}\,.\]
It follows by a direct calculation of the Laplace transform that, surprisingly,
\[ S_+\sim \Exp\{(1-\rho)(1-\lambda_+)\}\,.\]
In the case of $S_-$, we would have the arrival time
\[ A_i \sim \Exp'\{1-\rho\} + \ldots + \Exp'\{1-\rho\}\ \ \ \ \ ({\rm i.i.d. }\, k_-\ {\rm times})\,,\]
and the server time $B_i\sim \Exp\{\rho\}$. This means that for
\[ 1-\rho < \frac{k_-}{1+k_-}\,,\]
(note that $k_-/(1+k_-)$ is the fraction of particles to the left of zero), we define $\lambda_-=\lambda_-(\rho)$ as the smallest positive solution of
\begin{equation}\label{periodic-}
 \lambda_-(1 - \rho \lambda_-)^{k_-} = (1-\rho)^{k_-}\,,
 \end{equation}
and we get that
\[ S_- \sim \Exp\{\rho(1-\lambda_-)\}\,.\]
This easily implies that
\begin{equation}\label{periodic}
 \P(S_-\geq S_+) = \frac{(1-\rho)(1-\lambda_+)}{\rho(1-\lambda_-) + (1-\rho) (1-\lambda_+)}\,.
 \end{equation}
Therefore
$$\P(U\geq u)=\frac{(1-\rho_u)(1-\lambda_{+,u})}{\rho_u(1-\lambda_{-,u}) + (1-\rho_u) (1-\lambda_{+,u})}\,,$$
and
$$\P(\Theta\leq \alpha)=\frac{(1-\rho_\alpha)(1-\lambda_{+,\alpha})}{\rho_\alpha(1-\lambda_{-,\alpha}) + (1-\rho_\alpha) (1-\lambda_{+,\alpha})}\,,$$
(recall \eqref{periodic+} and \eqref{periodic-}). We remark that $k_+=k_-=\infty$ corresponds to $x=y=1$ in the previous example. In this case, $\lambda_+=\lambda_-=0$ and we get $(1-\rho)$ in the right-hand side of \eqref{periodic}, which agrees with the previous calculation.
\newline

\subsection{Bernoulli initial configuration} Choose $p_1\in[0,1)$ $p_2\in(p_1,1]$, and assume that $\{\eta(k)\,:\,k>0\}$ is a collection of i.i.d. Bernouilli random variables of parameter $p_1$, while $\{\eta(k)\,:\,k<0\}$, $\eta(k)$ is a collection of i.i.d. Bernouilli random variables of parameter $p_2$. In this case, we see that the profile $\sigma$ performs a two-sided random walk in $\ZZ^2$: when going down-right, $\sigma$ starts by moving down, and then it has probability $p_1$ to move down, and probability $1-p_1$ to move to the right. When going left-up, $\sigma$ starts by moving up, and then it has probability $p_2$ to move up, and probability $1-p_2$ to move to the left. We can again make the connection to a queue. Using the notation of Section \ref{sec:k+}, when considering $S^+$, we would get that $Z^+$ and $B_i$ (the first step and the server time) are distributed as a geometric sum of i.i.d. exponentials (the number of steps between two right-steps, so the expected number of down-steps is $1/(1-p_1)$), which is again exponential:
\[ Z_+, B_i \sim \Exp\{(1-p_1)(1-\rho)\}.\]
Furthermore, the arrival time $A_i$ is also a geometric sum of i.i.d. exponentials, so
\[ A_i \sim \Exp\{p_1\rho\}.\]
Note that in fact, in this case we have an M/M/1 queue. We need that $1-\rho>p_1$ for $S^+$ to be finite. It is not hard to repeat the argument of Section \ref{sec:k+} (other arguments for this classical problem are also possible): we find
\[ \lambda_+ = \frac{p_1}{1-p_1}\frac{\rho}{1-\rho}<1,\]
and from this
\[ S^+\sim \Exp\{1-\rho-p_1\}.\]
For $S_-$ we get the analogous result for $1-\rho<p_2$ by replacing $p_1$ by $1-p_2$ and $\rho$ by $1-\rho$:
\[ S^-\sim \Exp\{p_2-(1-\rho)\}.\]
This means that
\[ \P(S^-\geq S^+) = \frac{1-\rho-p_1}{p_2-p_1}.\]
Using that $\rho_u = (1+u)/2$, we find
$$\P\left(U\geq u\right)=\frac{(1-2p_1)-u}{2(p_2-p_1)}\,,$$
which is the uniform distribution on $(1-2p_2,1-2p_1)$, exactly as was found in \cite{FeKi}. The reader can also check that the above computation leads to the same limit distribution for the angle of the competition interface found in \cite{FeMaPi}.
\newline

\section{Computing the law of the asymptotic speed in the Hammersley process}\label{corollariesH}

\subsection{Periodic initial configuration} 
\begin{figure}[h]
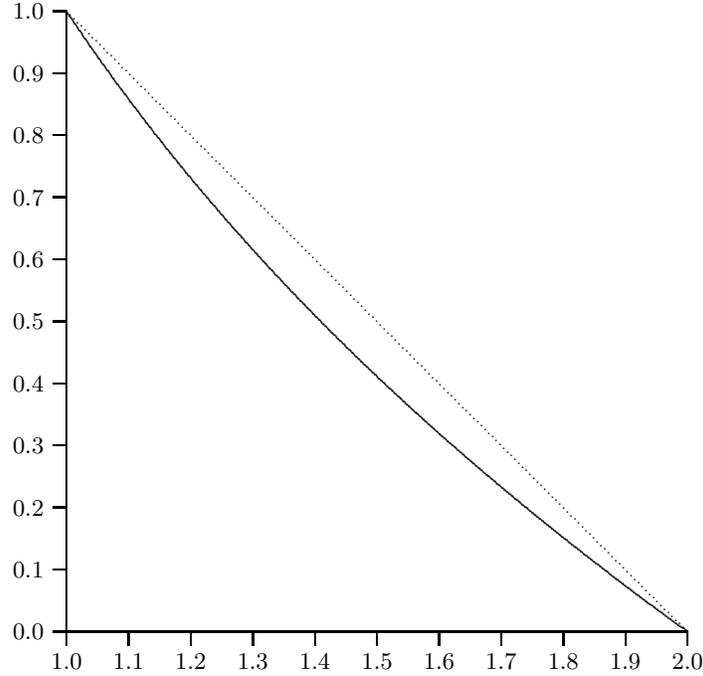

\[
\beginpicture
\footnotesize
\setcoordinatesystem units <0.5\textwidth,0.5\textwidth>
  \setplotarea x from 1 to 2, y from 0 to 1
  \axis bottom
    shiftedto y=0.0
    ticks numbered from 1 to 2 by 0.1
    /
  \axis left
    shiftedto x=1.0
    ticks numbered from 0 to 1 by 0.1
    /
\setquadratic
\plot
1.0000    1.0000     1.0050    0.9925     1.0100    0.9850    1.0150    0.9777    1.0200    0.9704    1.0250    0.9631    1.0300    0.9558    1.0350    0.9486    1.0400    0.9414    1.0450    0.9343    1.0500    0.9272    1.0550    0.9201    1.0600    0.9131    1.0650    0.9062    1.0700    0.8992    1.0750    0.8923    1.0800    0.8855    1.0850    0.8787    1.0900    0.8719    1.0950    0.8652    1.1000    0.8585    1.1050    0.8518    1.1100    0.8451    1.1150    0.8385
    1.1200    0.8320    1.1250    0.8255    1.1300    0.8190    1.1350    0.8125    1.1400    0.8061    1.1450    0.7997
    1.1500    0.7933    1.1550    0.7870    1.1600    0.7807    1.1650    0.7744    1.1700    0.7682    1.1750    0.7620
    1.1800    0.7558    1.1850    0.7496    1.1900    0.7435    1.1950    0.7374    1.2000    0.7314    1.2050    0.7253
    1.2100    0.7193    1.2150    0.7134    1.2200    0.7074    1.2250    0.7015    1.2300    0.6956    1.2350    0.6898
    1.2400    0.6839    1.2450    0.6781    1.2500    0.6723    1.2550    0.6666    1.2600    0.6608    1.2650    0.6551
    1.2700    0.6494    1.2750    0.6438    1.2800    0.6382    1.2850    0.6325    1.2900    0.6270    1.2950    0.6214
    1.3000    0.6159    1.3050    0.6104    1.3100    0.6049    1.3150    0.5994    1.3200    0.5940    1.3250    0.5885
    1.3300    0.5831    1.3350    0.5778    1.3400    0.5724    1.3450    0.5671    1.3500    0.5618    1.3550    0.5565
    1.3600    0.5512    1.3650    0.5460    1.3700    0.5407    1.3750    0.5355    1.3800    0.5303    1.3850    0.5252
    1.3900    0.5200    1.3950    0.5149    1.4000    0.5098    1.4050    0.5047    1.4100    0.4996    1.4150    0.4946
    1.4200    0.4895    1.4250    0.4845    1.4300    0.4795    1.4350    0.4746    1.4400    0.4696    1.4450    0.4647
    1.4500    0.4597    1.4550    0.4548    1.4600    0.4500    1.4650    0.4451    1.4700    0.4402    1.4750    0.4354
    1.4800    0.4306    1.4850    0.4258    1.4900    0.4210    1.4950    0.4162    1.5000    0.4115    1.5050    0.4067
    1.5100    0.4020    1.5150    0.3973    1.5200    0.3926    1.5250    0.3880    1.5300    0.3833    1.5350    0.3787
    1.5400    0.3740    1.5450    0.3694    1.5500    0.3648    1.5550    0.3602    1.5600    0.3557    1.5650    0.3511
    1.5700    0.3466    1.5750    0.3421    1.5800    0.3376    1.5850    0.3331    1.5900    0.3286    1.5950    0.3241
    1.6000    0.3197    1.6050    0.3152    1.6100    0.3108    1.6150    0.3064    1.6200    0.3020    1.6250    0.2976
    1.6300    0.2932    1.6350    0.2889    1.6400    0.2845    1.6450    0.2802    1.6500    0.2759    1.6550    0.2716
    1.6600    0.2673    1.6650    0.2630    1.6700    0.2587    1.6750    0.2545    1.6800    0.2502    1.6850    0.2460
    1.6900    0.2418    1.6950    0.2376    1.7000    0.2334    1.7050    0.2292    1.7100    0.2250    1.7150    0.2208
    1.7200    0.2167    1.7250    0.2125    1.7300    0.2084    1.7350    0.2043    1.7400    0.2002    1.7450    0.1961
    1.7500    0.1920    1.7550    0.1879    1.7600    0.1839    1.7650    0.1798    1.7700    0.1758    1.7750    0.1717
    1.7800    0.1677    1.7850    0.1637    1.7900    0.1597    1.7950    0.1557    1.8000    0.1517    1.8050    0.1478
    1.8100    0.1438    1.8150    0.1398    1.8200    0.1359    1.8250    0.1320    1.8300    0.1280    1.8350    0.1241
    1.8400    0.1202    1.8450    0.1163    1.8500    0.1125    1.8550    0.1086    1.8600    0.1047    1.8650    0.1009
    1.8700    0.0970    1.8750    0.0932    1.8800    0.0893    1.8850    0.0855    1.8900    0.0817    1.8950    0.0779
    1.9000    0.0741    1.9050    0.0703    1.9100    0.0666    1.9150    0.0628    1.9200    0.0590    1.9250    0.0553
    1.9300    0.0515    1.9350    0.0478    1.9400    0.0441    1.9450    0.0404    1.9500    0.0367    1.9550    0.0330
    1.9600    0.0293    1.9650    0.0256    1.9700    0.0219    1.9750    0.0182    1.9800    0.0146    1.9850    0.0109
    1.9900    0.0073    1.9950    0.0040    2.0000         0
/
\setdots <2pt>
\setlinear
\plot
1 1 2 0
/
\endpicture
\]
\caption{Picture of $P(S^+\geq S^-)$, together with uniform distribution.}\label{fig:detlamu}
\end{figure}

Choose $\lambda>0$ and $\mu>\lambda$, and assume that $(\nu^+(y)\,;\,y\geq 0)$ and $(\nu^-(y)\,;\,y<0)$ are deterministic periodic configurations of intensity $\lambda$ and $\mu$, i.e., $\nu^+$ is concentrated on $\{k/\lambda\,:\,k\geq 1\}$ while $\nu^-$ is concentrated on $\{k/\mu\,:\,k\leq -1\}$. Suppose $\lambda<\rho<\mu$. Theorem \ref{thm:2ndclasspart} then tells us that we should consider
\[ \P\left(\sup_{z\geq 0}  \{\nu(z) -\bar\nu_{\rho}(z)\} \geq \sup_{z<0} \{\nu(z) -\bar\nu_{\rho}(z)\}\right)\,,\]
where $\bar\nu_\rho$ is a Poisson process on $\R$. Define
\[ S^+ = \sup_{z\geq 0}  \{\nu(z) -\bar\nu_{\rho}(z)\}\ \ \ \mbox{and}\ \ \ S^- = \sup_{z<0} \{\nu(z) -\bar\nu_{\rho}(z)\}.\]
Clearly, for $z\geq 0$, we have $\nu(z) = \lfloor\lambda z\rfloor$, where $ \lfloor x\rfloor$ denotes the integer part of $x$. Therefore,
\begin{eqnarray*}
S^+ & = & \sup_{z\geq 0} \{ -\bar\nu_\rho(z) + \lfloor\lambda z\rfloor\}\\
& = & \lfloor\sup_{z\geq 0} \{ -\bar\nu_\rho(z) + \lambda z\}\rfloor\\
& = & -\lceil\inf_{z\geq 0} \{ \bar\nu_\rho(z) - \lambda z\}\rceil.
\end{eqnarray*}
This infimum of a Poisson process minus a linear function is studied in \cite{Py}. Theorem 3 in \cite{Py} entails that for $k\geq 0$
\begin{eqnarray*}
\P(S^+\geq k) & = & \P(\lceil\inf_{z\geq 0} \{ \bar\nu_\rho(z) - \lambda z\}\rceil \leq -k)\\
& = & \P(\inf_{z\geq 0} \{ \bar\nu_\rho(z) - \lambda z\}\leq -k)\\
& = & (1-p_+)^k,
\end{eqnarray*}
where $p_+=p_+(\lambda,\rho)\in (0,1)$ is the positive solution of
\[ p_+ = 1-e^{-p_+\rho/\lambda}.\]
For $z<0$ we have $\nu(z)=-\lfloor\mu |z|\rfloor$. So when we switch to positive $z$, we get
\begin{eqnarray*}
S^- & \stackrel{\cal D}{=} & \sup_{z\geq 0} \{\bar\nu_\rho(z) - \lfloor\mu z\rfloor\}\\
& = & \lceil\sup_{z\geq 0} \{\bar\nu_\rho(z) - \mu z\}\rceil.
\end{eqnarray*}
Now we can use Equation (7) of \cite{Py}: for $k\geq 0$
\begin{eqnarray*}
\P(S^-\leq k) & = & \P(\sup_{z\geq 0} \{\bar\nu_\rho(z) - \mu z\}\leq k)\\
& = & (1-\rho/\mu)\sum_{i=0}^k (-\rho/\mu)^i\frac{(k-i)^i}{i!}e^{\rho(k-i)/\mu}.
\end{eqnarray*}
In particular, $\P(S^- = 0) = 1-\rho/\mu$. Therefore,
\begin{eqnarray*}
\P(S^+\geq S^-) & = & \sum_{k=0}^\infty \P(S^-\leq k)p_+(1-p_+)^k\\
& = & p_+(1-\rho/\mu)\sum_{k=0}^\infty\sum_{i=0}^k (-\rho/\mu)^i\frac{(k-i)^i}{i!}e^{\rho(k-i)/\mu}(1-p_+)^k.
\end{eqnarray*}
We were not able to simplify this formula significantly. For $\lambda=1$ and $\mu=2$, Figure \ref{fig:detlamu} gives the graph of $P(S^+\geq S^-)$ as a function of $\rho\in (1,2)$.
We compare it with the probability $(\mu-\rho)/(\mu-\lambda)$, which we would get if we would take a Poisson process of intensity $\lambda$ for $x\geq 0$, and of intensity $\mu$ for $x\leq 0$ (see Section \ref{sec:Pois}).\\

\noindent An interesting limit is $\lambda\to\infty$. This means that $S^+=0$, and we get
\[\P(S^+\geq S_-) = 1-\rho/\mu.\]
This is exactly the same probability as when we would have a Poisson process of intensity $\mu$ left of $0$, and no particles to the right of $0$ (see Section \ref{sec:Pois})!

\subsection{Poisson initial configuration}\label{sec:Pois} Choose $\lambda>0$ and $\mu>\lambda$, and assume that $(\nu^+(y)\,;\,y\geq 0)$ and $(\nu^-(y)\,;\,y<0)$ are independent Poisson counting processes of intensity $\lambda$ and $\mu$, respectively. We will only consider the case $\rho\in(\lambda,\mu)$, the other cases are trivial. Define two asymmetric simple random walks $W^+$ and $W^-$, with
\[\P(W^+(n+1)-W^+(n) = +1)= p^+ := \frac{\lambda}{\lambda + \rho}\]
and
\[\P(W^-(n+1)-W^-(n) = +1)= p^- := \frac{\rho}{\mu + \rho}.\]
Since $\nu^+$, $\nu^+_{\rho}$, $\nu^-$ and $\nu^-_{\rho}$ are independent Poisson counting process, it is not hard to see that
$$\sup_{z\geq 0}\{ \nu^+(z)-\bar\nu^{+}_{\rho}(z)\} \stackrel{\D}{=} \sup_{n\geq 0} W^+(n)\,\,\mbox{ and that }\,\, \sup_{z< 0} \{\nu^-(z)-\bar\nu^{-}_{\rho}(z)\} \stackrel{\D}{=}  \sup_{n\geq 0} W^-(n).$$
Furthermore, it is well known (and easy to see) that $S^i:=\sup_{n\geq 0} W^i(n)\sim {\rm Geo}(r^i)$, where $r^i=p^i/(1-p^i)$, since $p^i<0.5$ (for $i=+,-$). We find
\begin{eqnarray*}
\P(S^+\geq S^-) & = & \sum_{k=0}^\infty \P(S^+\geq k)\P(S^-=k)\\
& = & \sum_{k=0}^\infty (r^+)^k(r^-)^{k}(1-r_-) =\frac{1-r^-}{1-r^+r^-} =\frac{\mu-\rho}{\mu-\lambda}\,.
\end{eqnarray*}
Therefore, by Theorem \ref{thm:2ndclasspart} (and taking $\rho_v=1\sqrt{v}$), this proves that
\[ \P(V^\nu \geq v) = \left\{\begin{array}{ll}
0 & \mbox{if }\,\, v\geq 1/ \lambda^2\\
\frac{\frac{1}{\sqrt{v}}-\lambda}{\mu-\lambda} & \mbox{if }\,\, 1/\mu^2< v< 1/\lambda^2\\
1 & \mbox{if }\,\,v\leq 1/\mu^2.
\end{array}
\right. \]
This agrees with the results found in \cite{CaDo,CoPi}.

\section{Final comments}
One interesting aspect of our method is that \eqref{comp-dist2} also holds in the lattice last-passage model with general i.i.d. weights, as soon as we can construct Busemann functions. In the exponential case, the limit shape is given by $f(u,v)=(\sqrt{u}+\sqrt{v})^2$. Also in the case of geometric weights, the shape function is explicitly known, and with our methods we could determine the distribution of the asymptotic speed of a second class particle. Busemann functions can be constructed under a curvature hypotheses for the respective limit shape \cite{Ne}. However, to prove that this curvature hypotheses holds in a great generality is by far one the most challenging problem in lattice last- (or first-) passage percolation models.

In the Hammersley last-passage model with random weights \cite{CaPi} (each Poissonian point $\x$ has a weight $X_\x$), the limit shape (due to the invariance of the model under volume preserving maps) $f(x,t)=\gamma\sqrt{xt}$ allows us to construct Busemann functions, and then prove \eqref{2class-dist} for the respective second class particle. However, we do not know the equilibrium distribution and so we can not go much further than that.

Part of the results concerning the distributional behavior of second class particles and competition interfaces in the rarefaction regime were already known  \cite{CaDo,CoPi,FeKi,FeMaPi,FePi,GuMo}. The genuine contributions  are \eqref{comp-dist1}, \eqref{comp-dist2} and \eqref{2class-dist}, and how it can be used to compute the distribution of the asymptotic speed, as soon as we have a good candidate for the equilibrium measure.
\newline

\paragraph{\bf Acknowledgments} We would like to thank Gerard Hooghiemstra for very useful discussions about the connections to queueing theory.

\end{document}